\documentclass[a4paper, twoside, notitlepage, 12pt]{amsart} 
\usepackage{amsmath,amsfonts,amssymb,amsthm,graphicx,mathrsfs} 

\usepackage{amscd}
\usepackage{xypic}
\xyoption{all}

\usepackage{verbatim}
\usepackage[latin1]{inputenc}   

\title[Symmetric Tensors and Divided Powers]{Counterexamples regarding Symmetric Tensors and Divided Powers}
\author{Christian Lundkvist} 
\address{Department of Mathematics \\ KTH \\ Stockholm \\ Sweden}
\email{chrislun@math.kth.se}

\let\ker\undefined
\DeclareMathOperator{\im}{Im}
\DeclareMathOperator{\ker}{Ker}

\DeclareMathOperator{\spec}{Spec}

\newcommand{\ra}{\longrightarrow}
\newcommand{\sra}{\rightarrow}

\newcommand{\NN}{\mathbb{N}}
\newcommand{\ZZ}{\mathbb{Z}}
\newcommand{\QQ}{\mathbb{Q}}

\newcommand{\pa}[1]{ \left( #1 \right) }
\newcommand{\mas}[1]{ \left\{ #1 \right\} }

\newcommand{\al}{\alpha}
\newcommand{\be}{\beta}
\newcommand{\ga}{\gamma}

\newcommand{\fii}{\varphi}
\newcommand{\sig}{\sigma}
\newcommand{\Ga}{\Gamma}

\newcommand{\eqbeg}{\begin{equation}}
\newcommand{\eqend}{\end{equation}}
\newcommand{\arbeg}[1]{\begin{array}{#1}}
\newcommand{\arend}{\end{array}}

\newcommand{\sse}{\subseteq}

\newcommand{\ten}{\otimes}

\newcommand{\tim}{\times}

\newcommand{\Sets}{\mathbf{Sets}}

\newcommand{\Alg}[1]{#1\textrm{-}\mathbf{Alg}}
\newcommand{\Mod}[1]{#1\textrm{-}\mathbf{Mod}}

\newcommand{\FF}{\mathcal{F}}

\renewcommand{\hom}{\mathrm{Hom}}

\newtheorem{lemma}[subsection]{Lemma}
\newtheorem{cor}[subsection]{Corollary}
\newtheorem{prop}[subsection]{Proposition}

\theoremstyle{definition}
\newtheorem{defi}[subsection]{Definition}
\newtheorem{exam}[subsection]{Example}

\newtheorem{para}[subsection]{}
\theoremstyle{remark}
\newtheorem{rem}[subsection]{Remark}

\numberwithin{equation}{subsection}

\newcommand{\TS}{\mathrm{TS}}
\newcommand{\T}{\mathrm{T}}
\newcommand{\rmS}{\mathrm{S}}

\newcommand{\calI}{\mathcal{I}}
\newcommand{\fS}{\mathfrak{S}}

\newcommand{\Pol}{\mathrm{Pol}}
\newcommand{\mdeg}{\mathrm{mdeg}}

\begin{document}
\pagenumbering{arabic}
\maketitle
\begin{abstract} 
We investigate the similarities and differences between the module of symmetric tensors $\TS^n_A(M)$ and the module of divided powers $\Ga^n_A(M)$. There is a canonical map $\Ga^n_A(M) \sra \TS^n_A(M)$ which is an isomorphism in many important cases. We give examples showing that this map need neither be surjective nor injective in general. These examples also show that the functor $\TS_A^n$ does not in general commute with base change.
\end{abstract}


\section*{Introduction}

Symmetric tensors and divided powers are important tools in algebraic geometry. They appear for instance in the study of Hilbert and Chow schemes parametrizing zero-dimensional subschemes or cycles of a given scheme (c.f. \cite{fga_hilb_schemes}, \cite[5.5]{sga4_deligne_coh_supp_prop}, \cite{rydh_representability_of_gamma_inprep}, \cite{rydh_gammasymchow_inprep}, \cite{rydh-skjelnes_gen_etale_fam}).

For a flat family of schemes the symmetric tensors and divided powers coincide. However, for non-flat families they may differ, and it is then of interest to understand the relations between the resulting schemes. 

The purpose of this article is to provide examples showing how symmetric tensors and divided powers may differ in the non-flat case. We shall throughout the article stick to the case of affine schemes.\\
\\
Let $A$ be a commutative ring with unit element. The module of symmetric tensors $\TS^n_A(M)$ for an $A$-module $M$ is defined as the submodule of elements of $M \ten_A \cdots \ten_A M$ invariant under the natural action of the symmetric group $\fS_n$. When $A$ is a field of characteristic zero these objects have been studied since the nineteenth century (see e.g. \cite{junker_93}).

More recently a related object has been introduced, the module of divided powers $\Ga^n_A(M)$ \cite{roby_lois_pol}. This module is not defined as intuitively as $\TS^n_A(M)$, but it is functorially more well-behaved. For instance $\Ga^n_A(M)$ satisfies a universal property regarding polynomial laws, and commutes with arbitrary base change $A \sra A'$. The module of symmetric tensors on the other hand commutes with flat base change $A \sra A'$ but not any base change in general. This has been pointed out in \cite[5.5.2.7]{sga4_deligne_coh_supp_prop} but the author does not know of any published counterexamples.

There is a canonical map $\Ga^n_A(M) \sra \TS^n_A(M)$ comparing the two modules. This map is an isomorphism when $n!$ is invertible in $A$, or when $M$ is a flat $A$-module. The purpose of this article is to give examples showing that the map $\Ga^n_A(M) \sra \TS^n_A(M)$ is in general neither injective nor surjective (Examples \ref{ex:rydh_canmap_not_injective}, \ref{ex:canmap_not_injective}, \ref{ex:canmap_not_surj_1}, \ref{ex:canmap_not_surj_2}). These examples also show that the functor $\TS^n_A$ does not commute with base change in general. Specifically, we show that the base change map $\TS^n_A(M) \ten_A A' \sra \TS^n_{A'}(M \ten_A A')$ is neither injective nor surjective in general (Examples \ref{ex:TS_base_change_not_injective}, \ref{ex:TS_base_change_not_surjective}).

Furthermore we show in Section \ref{sec:from_modules_to_algebras} that if the module $M$ has the property that the canonical map fails to be injective/surjective, then the symmetric algebra $\rmS_A(M)$ will also have this property. Thus the examples are extended from modules to graded algebras.

Finally, in Section \ref{sec:alg_struct_on_ext_pow} we relate our examples to work of Laksov and Thorup \cite{laksovthorup_determinantal} who discuss the structure of the exterior product $\bigwedge^n_A(M)$ as a module over $\TS^n_A(B)$, where $B$ is an $A$-algebra. This module structure then gives formulas related to Schubert calculus and intersection theory of flag schemes. We use Example \ref{ex:canmap_not_surj_2} to show that $\bigwedge^n_A(M)$ does not in general admit a structure of $\TS^n_A(B)$-module.

\textbf{Acknowledgements.} I would like to thank Francesco Vaccarino and David Rydh for helpful suggestions and remarks. 

\section{Definitions and first properties}

For the convenience of the reader we present some definitions and results concerning symmetric tensors and divided powers. All information in this section can be found in \cite{roby_lois_pol} or \cite{ferrand_norme}.

For the remainder of this section, fix a commutative ring $A$ with unit element, an $A$-module $M$ and an integer $n$. For the rest of the paper all rings will be assumed to be commutative with identity unless otherwise specified.

\begin{defi}[Symmetric tensors]
Denote by $\T^n_A (M)$ the $n$-fold tensor product 
\[
\T^n_A (M) = \underbrace{M \ten_A \ldots \ten_A M}_{n}. 
\]
The tensor product $\T^n_A(M)$ has a canonical $A$-module structure and the symmetric group $\fS_n$ acts on $\T^n_A (M)$ by $A$-module homomorphisms defined by 
\[
\sig(m_1 \ten \ldots \ten m_n) = m_{\sig^{-1}(1)} \ten \ldots \ten m_{\sig^{-1}(n)} 
\]
for $\sig \in \fS_n$. If $M$ is free with basis $\mas{e_i}_{i \in \calI}$, then $\T^n_A(M)$ is free with basis $\mas{e_{i_1} \ten \ldots \ten e_{i_n}}_{(i_1, \ldots, i_n) \in \calI^n}$.

The module of invariants $\T^n_A(M)^{\fS_n}$ is called the \emph{module of symmetric tensors} and is denoted $\TS^n_A(M)$. 
\end{defi}

\begin{defi}[Shuffle product]
\label{def:shuffle_prod}
Consider the direct sum $\TS_A(M) = \bigoplus_{k \geq 0} \TS^k_A(M)$. We have a product $\tim$ on $\TS_A(M)$ called the \emph{shuffle product} which makes $\TS_A(M)$ into a commutative graded ring. The product is defined as follows: Let $z \in \TS^k_A(M)$ and $z' \in \TS^l_A(M)$. Then
\[
z \tim z' = \sum_{\sig \in \fS_{k,l}} \sig(z \ten z')
\]
where $\fS_{k,l}$ is the subset of elements $\sig \in \fS_{k+l}$ such that $\sig(1) < \sig(2) < \ldots < \sig(k)$ and $\sig(k+1) < \sig(k+2) < \ldots < \sig(k+l)$.
\end{defi}

\begin{defi}[Polynomial laws]
Let $N$ be an $A$-module. A \emph{polynomial law} from $M$ to $N$ is defined as follows: Let $\FF_M: \Alg{A} \sra \Sets$ be the functor defined by $\FF_M(A') = M \ten_A A'$ viewed as a \emph{set}. Then a polynomial law $F$ from $M$ to $N$ is a natural transformation $\FF_M \sra \FF_N$.

In other words, for each morphism of $A$-algebras $g:A' \sra A''$ we have a commutative diagram
\[
\xymatrix{
M \ten_A A' \ar[r]^{F_{A'}} \ar[d]_{1_M \ten g} &  N \ten_A A' \ar[d]^{1_N \ten g}\\
M \ten_A A'' \ar[r]_{F_{A''}} &  N \ten_A A''\\
}
\]
where the horizontal maps are maps of the underlying sets, and not homomorphisms of modules in general.

The polynomial law $F$ is called \emph{homogeneous of degree} $n$ if $F_{A'}(ax) = a^n F_{A'}(x)$ for each $a \in A'$ and each $x \in M \ten_A A'$.

If $B$ and $C$ are (not necessarily commutative) $A$-algebras, then a polynomial law $F: B \sra C$ is called \emph{multiplicative} if $F_{A'}(xy) = F_{A'}(x) F_{A'}(y)$ for each $x,y \in B \ten_A A'$.
\end{defi}

\begin{defi}[Divided powers]
For an $A$-module $M$ there exists a commutative graded algebra $\Ga_A(M) = \bigoplus_{n \geq 0} \Ga^n_A(M)$ with multiplication $\tim$, together with set maps $\ga^n: M \sra \Ga^n_A(M)$ such that for each $a \in A$, $x,y \in M$ and $n,m \in \NN$ we have
\[
\arbeg{rcl}
\Ga^0_A(M) & = & A \ \ \ \textrm{and} \ \ \ga^0(x) = 1,\\
\Ga^1_A(M) & = & M \ \ \ \textrm{and} \ \ \ga^1(x) = x,\\
\ga^n(ax) & = & a^n \ga^n(x),\\
\ga^n(x+y) & = & \sum_{i=0}^n \ga^{i}(x) \tim \ga^{n-i}(y),\\
\ga^n(x) \tim \ga^m(x) & = & {n+m \choose n} \ga^{n+m}(x).
\arend
\]
If $(x_i)_{i \in \calI}$ is a family of elements of $M$, and $\nu = (\nu_i)_{i \in \calI}$ is a multiindex of finite support, then we write
\[
\ga^{\nu}(x) := \mathop{\tim}_{i \in \calI} \ga^{\nu_i}(x_i). 
\]
We have that $\ga^{\nu}(x) \in \Ga^n_A(M)$ where $n = |\nu| = \sum_{i \in \calI} \nu_i$.
\end{defi}

\begin{para}[Functoriality]
The application $M \mapsto \Ga_A(M)$ is a functor from $A$-modules to graded $A$-algebras \cite[Ch. III \S 4, p. 251]{roby_lois_pol}. 
\end{para}

\begin{para}[Base change]
\label{par:Ga_base_change}
For each morphism $A \sra A'$ there is a natural map
\[
\Ga_A(M) \ten_A A' \ra  \Ga_{A'}(M \ten_A A') 
\]
defined by $\ga^n(x) \ten 1 \mapsto \ga^n(x \ten 1)$, which is an isomorphism \cite[Thm. III.3, p. 262]{roby_lois_pol}. Thus the maps $\ga^n_{A'}: M \ten_A A' \sra \Ga^n_{A'}(M \ten_A A')$ define a polynomial law $\ga^n: M \sra \Ga^n_A(M)$. This polynomial law is homogeneous of degree $n$.
\end{para}

\begin{para}[Universal property]
\label{par:univ_prop}
For $A$-modules $M$, $N$ we write $\Pol^n_A(M,N)$ for the set of polynomial laws $M \sra N$ of degree $n$. Then the natural map $\hom_A(\Ga^n_A(M),N) \sra \Pol^n_A(M,N)$ given by $f \mapsto f \circ \ga^n$ is an isomorphism. Thus $\Ga^n_A(M)$ represents the functor $N \mapsto \Pol^n_A(M,N)$.
\end{para}

\section{The canonical map}

In this section we define the canonical map $\Ga^n_A(M) \sra \TS^n_A(M)$ and give critera for when this map is injective or surjective.

\begin{defi}
\label{def:Ga_TS_canmap}
Let $A$ be a ring, $n$ an integer and $M$ an $A$-module. There is a homogeneous polynomial law of degree $n$ from $M$ to $\TS^n_A(M)$ defined by sending an element $x \in M$ to $x^{\ten n} \in \TS^n_A(M)$.

By (\ref{par:univ_prop}) this polynomial law gives rise to an $A$-module homomorphism 
\[
\Ga^n_A(M) \sra \TS^n_A(M) 
\]
that maps $\ga^n(x)$ to $x^{\ten n}$ for $x \in M$.
\end{defi}

\begin{prop}
\label{prp:Ga_TS_iso}
The morphism of Definition \ref{def:Ga_TS_canmap} is an isomorphism in the following important cases:
\begin{itemize}
\item[(i)] The element $n!$ is invertible in the ring $A$ \cite[Prop. III.3, p. 256]{roby_lois_pol}.
\item[(ii)] The $A$-module $M$ is free \cite[Prop. IV.5, p. 272]{roby_lois_pol}.
\item[(iii)] More generally, when the $A$-module $M$ is flat \cite[5.5.2.5, p. 123]{sga4_deligne_coh_supp_prop}.
\end{itemize}
\end{prop}

\begin{para}[Factorization of the canonical morphism]
\label{par:fact_of_can_mor}
Let $A$ be a ring and $M$ an $A$-module with presentation
\[
0 \ra P \ra F \ra M \ra 0
\]
with $F$ a free $A$-module. Then the surjection $F \sra M$ induces a surjection $\Ga^n_A(F) \sra \Ga^n_A(M)$ with kernel $K$ given by
\[
K = \langle \ga^s(p) \tim y: \ p \in P, \ y \in \Ga^{n-s}_A(F), \ 1 \leq s \leq n \rangle,
\]
by \cite[Prop. IV.8, p. 284]{roby_lois_pol}. Since $\Ga^n_A(F) \cong \TS^n_A(F)$ by Proposition~\ref{prp:Ga_TS_iso} we can view $K$ as a submodule of $\TS^n_A(F)$, and with this identification we then have $\Ga^n_A(M) \cong \TS^n_A(F)/K$.

Let $N \sse \T^n_A(F)$ denote the kernel of the map $\pi: \T^n_A(F) \sra \T^n_A(M)$. Then $N$ is stable under the action of $\fS_n$. Furthermore, the functor $(\cdot)^{\fS_n}$ is left exact so the exact sequence
\[
0 \ra N \ra \T^n_A(F) \ra \T^n_A(M) \ra 0
\]
gives an exact sequence
\[
0 \ra N^{\fS_n} \ra \TS^n_A(F) \ra \TS^n_A(M).
\]
Thus we have a canonical injection $\TS^n_A(F) / N^{\fS_n} \sra \TS^n_A(M)$.

Also, we note that $K \sse N^{\fS_n}$, and so we have a surjection $\TS^n_A(F)/K \sra \TS^n_A(F) / N^{\fS_n}$. Thus, the canonical map $\Ga^n_A(M) \sra \TS^n_A(M)$ factors as
\eqbeg
\label{eq:fact_of_can_map}
\Ga^n_A(M) \cong \TS^n_A(F)/K \ra \TS^n_A(F) / N^{\fS_n} \ra \TS^n_A(M)
\eqend
where the first map is surjective and the second is injective.
\end{para}

\begin{prop}
\label{prp:canmap_surj_inj}
With the notation of (\ref{par:fact_of_can_mor}), we have that the map $\Ga^n_A(M) \sra \TS^n_A(M)$ is
\begin{itemize}
\item[$(a)$] injective if and only if $K = N^{\fS_n}$,
\item[$(b)$] surjective if and only if $\TS^n_A(F) \sra \TS^n_A(M)$ is surjective. Moreover, the image of $\TS^n_A(F)$ in $\TS^n_A(M)$ is generated by the elements
\[
m_{\nu} := \mathop{\tim}_{i \in \calI} m_i^{\ten \nu_i}
\]
where $\mas{m_i}_{i \in \calI}$ is any prescribed generating set of $M$, $\nu$ is a multiindex of finite support and $|\nu| = \sum_{i \in \calI} \nu_i = n$. Here $\tim$ denotes the shuffle product of Definition \ref{def:shuffle_prod}.
\end{itemize}
\end{prop}

\begin{proof}
To prove $(a)$, we note that by the factorization (\ref{eq:fact_of_can_map}) we have that $\Ga^n_A(M) \sra \TS^n_A(M)$ is injective if and only if $\TS^n_A(F)/K \sra \TS^n_A(F)/ N^{\fS_n}$ is an isomorphism. This happens if and only if $K = N^{\fS_n}$.\\
\\
For $(b)$ we have that the factorization (\ref{eq:fact_of_can_map}) further implies that $\Ga^n_A(M) \sra \TS^n_A(M)$ is surjective if and only if $\TS^n_A(F)/N^{\fS_n} \sra \TS^n_A(M)$ is an isomorphism. This happens if and only if $\TS^n_A(F) \sra \TS^n_A(M)$ is surjective.

To show the last part of $(b)$, suppose that $\mas{e_i}_{i \in \calI}$ is a basis for $F$ and that $F \sra M$ maps $e_i$ to $m_i$ for all $i \in \calI$. Then the corresponding elements $e_{\nu} := \tim_{i \in \calI} e_i^{\ten \nu_i}$ with $|\nu|=n$ form a basis for $\TS^n_A(F)$ \cite[IV §5 Prop. 4]{bourbaki_algebre_4_5}. The images of the elements $e_{\nu} \in \TS^n_A(F)$ are the elements $m_{\nu} \in \TS^n_A(M)$.
\end{proof}

\section{Injectivity of the canonical map}

Here we give two examples showing that the map $\Ga^n_A(M) \sra \TS^n_A(M)$ need not be injective.

\begin{exam}
\label{ex:rydh_canmap_not_injective}
This short example of non-injectivity is due to David Rydh \cite{rydh_gammasymchow_inprep}. Recall that if $p$ is a prime, then $p \mid {p \choose s}$ whenever $1 \leq s < p$.

Let $k$ be a field of prime characteristic $p$, and let $A = k[x]$ and $B = k$, where $A \sra B$ sends $x$ to $0$. Then $\T^p_A(B) \cong k$ and so $\TS^p_A(B) \cong k$. However, we have $\Ga^p_A(B) \cong \TS^p_A(A)/K$ by (\ref{par:fact_of_can_mor}) where 
\[
K = \langle x^{\ten s} \tim 1^{\ten(p-s)}: \ 1 \leq s \leq p \rangle = \langle x^s{p \choose s} 1^{\ten p}: \ 1 \leq s \leq p\rangle = \langle x^p 1^{\ten p} \rangle. 
\]
Thus $\Ga^p_A(B) \cong k[x]/(x^p)$ and hence the map $\Ga^p_A(B) \sra \TS^p_A(B)$ is not injective.
\end{exam}

\begin{exam}
\label{ex:canmap_not_injective}
This example gives a morphism of rings $A \sra A'$ and an $A$-module $M$ such that $\Ga^n_A(M) \sra \TS^n_A(M)$ is injective, while $\Ga^n_{A'}(M') \sra \TS^n_{A'}(M')$ is not injective, where $M' = M \ten_A A'$.

Let $k$ be a field of characteristic $2$, and let $A = k[s,t]$ be the polynomial ring in two variables $s,t$. Moreover, let $A'$ be the algebra $A' = k[s,t,z]/(z(s+t))$. Consider the free module $F = A^2$ with generators $e_1, e_2$ and let $M = F/ \langle n \rangle$, where $n = s e_1 + t e_2$. Let $m_1, m_2$ be the images of $e_1, e_2$ in $M$ and denote by $M'$ the module $M \ten_A A'$.\\
\\
First we show that $\Ga^2_A(M) \sra \TS^2_A(M)$ is injective. By Proposition \ref{prp:canmap_surj_inj} we thus need to check that $K = N^{\fS_2}$, where $K$ is the kernel of 
\[
\TS^2_A(F) \cong \Ga^2_A(F) \sra \Ga^2_A(M)
\]
and $N$ is the kernel of the map $\T^2_A(F) \sra \T^2_A(M)$. By (\ref{par:fact_of_can_mor}) we have that
\[
K = \langle n \tim e_1, \ n \tim e_2, \ n^{\ten 2} \rangle = \langle n \ten e_1 + e_1 \ten n, \  n \ten e_2 + e_2 \ten n, \ n \ten n \rangle.
\]

To compute $N^{\fS_2}$ we first note that $N$ is generated by the elements
\[
\mas{e_1 \ten n, \ n \ten e_1, \ e_2 \ten n, \ n \ten e_2}.
\]
Choose an element $u \in N^{\fS_2} = N \cap \TS^2(F)$ and let $\sig: \T_A^2(F) \sra \T_A^2(F)$ be the homomorphism defined by $\sig(e_i \ten e_j) = e_j \ten e_i$ for $i,j = 1,2$. We write $u$ as
\[
u = a n \ten e_1 + b e_1 \ten n + c n \ten e_2 + d e_2 \ten n,
\]
where $a,b,c,d \in A = k[s,t]$. We have $u + \sig(u) = u + u = 0$, and so
\[
0 = (a+b)(n \ten e_1 + e_1 \ten n) + (c+d)(n \ten e_2 + e_2 \ten n).
\]
Using that $n = s e_1 + t e_2$ and cancelling terms we obtain
\[
0 = ((a+b)t + (c+d)s)(e_1 \ten e_2 + e_2 \ten e_1).
\]
Hence
\eqbeg
\label{eq:kerTSF_TSM2}
(a+b)t + (c+d)s = 0
\eqend
and we conclude that $s | (a+b)$. Hence $a+b = fs$ for some $f \in A$, and from (\ref{eq:kerTSF_TSM2}) we obtain $(c+d)s = fts$ and so $c+d = ft$. We conclude that $b = a + fs$ and $d = c+ft$ and so $u$ can be written as
\[
u = a n \ten e_1 + (a+fs) e_1 \ten n + c n \ten e_2 + (c+ft) e_2 \ten n = 
\]
\[
= a (n \ten e_1 + e_1 \ten n) + c(n \ten e_2 + e_2 \ten n) + f n \ten n.
\]
Thus $N^{\fS_2}$ is generated by the elements
\[
\mas{(n \ten e_1 + e_1 \ten n), \ (n \ten e_2 + e_2 \ten n), \ n \ten n }
\]
and so $K = N^{\fS_2}$. Hence $\Ga^2_A(M) \sra \TS^2_A(M)$ is injective.\\ 
\\
Next we show that $\Ga^2_{A'}(M') \sra \TS^2_{A'}(M')$ is not injective. Let $K'$ denote the kernel of 
\[
\TS^2_{A'}(F') \cong \Ga^2_{A'}(F') \sra \Ga^2_{A'}(M')
\]
and denote by $N'$ the kernel of $\T^2_{A'}(F') \sra \T^2_{A'}(M')$. We will show that $K' \subset (N')^{\fS_2}$ is a proper subset.

The element
\[
v = zs(e_1 \ten e_1 + e_2 \ten e_2) = zt(e_1 \ten e_1 + e_2 \ten e_2) \in \T^2_{A'}(F')
\]
is clearly in $\TS^2_{A'}(F')$. In $M'$ we have $zsm_1 = ztm_2 = zsm_2$ and so the image of $v$ under the map $\T^2_{A'}(F') \sra \T^2_{A'}(M')$ is
\[
zs(m_1 \ten m_1 + m_2 \ten m_2) = zs m_2 \ten m_2 + zs m_2 \ten m_2 = 0.
\]
Thus $v \in (N')^{\fS_2}$. Assume, to obtain a contradiction, that $v \in K'$. We have by (\ref{par:fact_of_can_mor}) that
\[
K' = \langle n \ten e_1 + e_1 \ten n, \  n \ten e_2 + e_2 \ten n, \ n \ten n \rangle
\]
and we see that we can choose generators as
\[
K' = \langle t(e_2 \ten e_1 + e_1 \ten e_2), \ s(e_1 \ten e_2 + e_2 \ten e_1), \ s^2 e_1 \ten e_1 + t^2 e_2 \ten e_2 \rangle. 
\]
We have by \cite[IV §5 Prop. 4]{bourbaki_algebre_4_5} that $\TS^2_{A'}(F')$ is a free $A'$-module of rank $3$ generated by the elements
\[
f_1 = e_1 \ten e_1, \ \ f_2 = e_2 \ten e_2, \ \ f_{12} = e_1 \ten e_2 + e_2 \ten e_1.
\]
With this notation we have
\[
K' = \langle s f_{12}, \ t f_{12}, \ s^2 f_1 + t^2 f_2 \rangle.
\]
Now let $B = k[s,t,z]$ be the polynomial ring and let $G = B^3$ be a free module with basis $f_1, f_2, f_{12}$. Then we have $v \in K'$ if and only if
\[
zs(f_1 + f_2) \in \langle s f_{12}, \ t f_{12}, \ s^2 f_1 + t^2 f_2, \ z(s+t)f_1, \ z(s+t)f_2, \ z(s+t)f_{12} \rangle
\]
where the above are elements of the free $B$-module $G$.
 
Thus
\[
zs(f_1 + f_2) = a (s^2 f_1 + t^2 f_2) + b z(s+t) f_1 + c z(s+t) f_2 =
\]
\[
= (as^2 + bz(s+t))f_1 + (at^2 + cz(s+t))f_2,
\]
where $a,b,c \in B$. Comparing terms on each side, we conclude that 
\eqbeg
\label{eq:final_step_counterexample_injectivity}
zs = as^2 + bz(s+t).
\eqend
From this we have $z \mid as^2$ and so $z \mid a$. By the same reason we have that $s \mid b$. Hence the polynomial on the right hand side of (\ref{eq:final_step_counterexample_injectivity}) is either zero or has degree $\geq 3$, a contradiction. We conclude that $v \notin K'$, and thus the inclusion $K' \subset (N')^{\fS_2}$ is strict. Hence $\Ga^2_{A'}(M') \sra \TS^2_{A'}(M')$ is not injective.
\end{exam}

\begin{rem}
\label{rem:canmap_not_injective}
It is possible to extend the non-injectivity part of Example \ref{ex:canmap_not_injective} to characteristic $p \geq 2$ as follows: Let $k$ be a field of characteristic $p$ and let $A = k[s_1, \ldots, s_p]$ be the polynomial ring in $p$ variables. Define
\[
A' = k[s_1, \ldots, s_p, z]/(zs_1 - zs_i: \ 2 \leq i \leq p),
\]
with the obvious map $A \to A'$. Let $F = A^p$ with basis $e_1, \ldots, e_p$ and let 
\[
M = F / \langle s_1 e_1 - s_i e_i: \ 2 \leq i \leq p \rangle.
\]
The goal is now to show that $\Ga^p_A(M) \to \TS^p_A(M)$ is injective while $\Ga^p_{A'}(M') \to \TS^p_{A'}(M')$ is not injective, where $M' = M \ten_A A'$. The non-injectivity of $\Ga^p_{A'}(M') \to \TS^p_{A'}(M')$ is shown as follows: Denote by $N'$ the kernel of $\T^p_{A'}(F') \to \T^p_{A'}(M')$. Then the element
\[
v = z s_1 (e_1^{\ten p} + \ldots + e_p^{\ten p}) \in \T^p_{A'}(F')
\]
is in $(N')^{\fS_p}$. However, the elements of $K' = \ker(\TS^p_{A'}(F') \to \TS^p_{A'}(M'))$ containing terms of the form $a_i e_i^{\ten p}$ must satisfy $a_i = s_i^p b_i$ with $b_i \in A'$. Thus $v \notin K'$ by reasons of homogeneity for $p \geq 3$, while the case $p = 2$ is already in the Example. This shows that $K' \ne (N')^{\fS_p}$ and so we have shown non-injectivity.

The map $\Ga^p_A(M) \to \TS^p_A(M)$ is probably injective, but it is unclear how to extend the methods of Example \ref{ex:canmap_not_injective} to show this.
\end{rem}

\section{Surjectivity of the canonical map}

In this section we give two lemmas which give special cases where the map $\Ga^n_A(M) \sra \TS^n_A(M)$ is surjective. We also give an algorithmic method of checking surjectivity, and finally we provide two examples showing that the canonical map need not be surjective in general.

\begin{lemma}
\label{lem:2gensTS2}
Let $M$ be an $A$-module generated by two elements $m_1, m_2$. Then the morphism
\[
\Ga^2_A(M) \ra \TS^2_A(M)
\]
is surjective.
\end{lemma}

\begin{proof}
By Proposition \ref{prp:canmap_surj_inj} it is enough to show that $\TS^2_A(M)$ is generated by the elements
\eqbeg
\label{eq:3elements}
m_1 \ten m_1, \quad m_2 \ten m_2, \quad m_1 \ten m_2 + m_2 \ten m_1.
\eqend
Let $u \in \TS^2_A(M)$ be any element. This element can be written as
\[
u = a_{11} m_1 \ten m_1 + a_{22} m_2 \ten m_2 + a_{12} m_1 \ten m_2 + a_{21} m_2 \ten m_1
\]
with $a_{ij} \in A$. We write
\[
u = a_{11} m_1 \ten m_1 + a_{22} m_2 \ten m_2 + a_{21} (m_1 \ten m_2 + m_2 \ten m_1) + 
(a_{12}-a_{21}) m_1 \ten m_2.
\]
It is clear from the above that the element $(a_{12}-a_{21})m_1 \ten m_2$ is in $\TS^2_A(M)$, so we are done if we show that this element is a linear combination of the three elements (\ref{eq:3elements}).

Let $a = a_{12}-a_{21}$, and denote by $F$ the free module $A^2$ generated by the basis elements $e_1, e_2$. Then $M$ is isomorphic to a quotient $F/N$ where $N \sse F$ is generated by elements $\mas{f_i e_1 - g_i e_2}_{i \in \calI}$ with $f_i, g_i \in A$. The isomorphism is given by $e_i \mapsto m_i$ for $i=1,2$.

Let $n_i = f_i e_1 - g_i e_2$. Then $M \ten_A M \cong (F \ten_A F) / N'$, where $N'$ is the module generated by the elements
\[
\mas{n_i \ten e_1, n_i \ten e_2, e_1 \ten n_i, e_2 \ten n_i}_{i \in \calI}.
\]
Since the element $a m_1 \ten m_2 - a m_2 \ten m_1$ is zero in $M \ten_A M$ we conclude that 
\[
a e_1 \ten e_2 - a e_2 \ten e_1 \in N'
\]
and we therefore have
\eqbeg
\label{eq:m1_ten_m2}
a e_1 \ten e_2 - a e_2 \ten e_1 = 
\sum_{i \in \calI} (x_i n_i \ten e_1 + y_i n_i \ten e_2 + z_i e_1 \ten n_i + w_i e_2 \ten n_i)
\eqend
where the elements $x_i,y_i,z_i,w_i$ are in $A$ and only a finite number of these elements are non-zero. Inserting $n_i = f_i e_1 - g_i e_2$ in (\ref{eq:m1_ten_m2}) and comparing the coefficients of $e_1 \ten e_2$ we obtain
\[
a = \sum_{i \in \calI} (y_i f_i - z_i g_i).
\]
Since $f_i m_1 = g_i m_2$ in $M$ we have
\[
a m_1 \ten m_2 = \sum_{i \in \calI} y_i f_i m_1 \ten m_2 - \sum_{i \in \calI} z_i g_i m_1 \ten m_2 = 
\]
\[
= \sum_{i \in \calI} y_i g_i m_2 \ten m_2 - \sum_{i \in \calI} z_i f_i m_1 \ten m_1.
\]
This shows that $a m_1 \ten m_2$ is a linear combination of the elements (\ref{eq:3elements}).
\end{proof}

\begin{lemma}
\label{lem:TS2_UFD}
Let $A$ be a UFD and let $M$ be a module of the form $A^k/\langle f \rangle$, where $f$ is defined as $f = \sum_{i=1}^k f_i e_i \in A^k$ with $f_i \in A$ and $\mas{e_i}$ is the canonical basis of $A^k$. Suppose further that $\gcd(f_k,f_i) = 1$ for $i = 1, \ldots, k-1$. Then
\[
\Ga^2_A(M) \sra \TS^2_A(M)
\]
is surjective.
\end{lemma}

\begin{proof}
Let $F = A^k$, and let $F \sra M$ be the canonical surjective map sending $e_i$ to $m_i$ where $\mas{m_i}$ is a set of generators of $M$. By Proposition \ref{prp:canmap_surj_inj} we need to check that $\TS^2_A(M)$ is generated by the elements
\eqbeg
\label{eq:TS2_UFD_generators1}
m_i \ten m_j + m_j \ten m_i, \quad m_i \ten m_i, \quad 1 \leq i,j \leq k.
\eqend
First we wish to show that the elements $\mas{m_i \ten m_j}$ with $1 \leq i,j \leq k-1$ are linearly independent. This linear independence implies that the submodule $L \sse \T^2_A(M)$ generated by $\mas{m_i \ten m_j}_{i,j \leq k-1}$ is isomorphic to $\T^2_A(A^{k-1})$, and hence that the module of invariants $L^{\fS_2} \sse \TS^2_A(M)$ is isomorphic to $\TS^2_A(A^{k-1})$. Thus by Proposition \ref{prp:canmap_surj_inj} and Proposition \ref{prp:Ga_TS_iso} the elements of $L^{\fS_2}$ can be generated by 
\eqbeg
\label{eq:TS2_UFD_generators2}
m_i \ten m_j + m_j \ten m_i, \quad m_i \ten m_i, \quad 1 \leq i,j \leq k-1.
\eqend
Let $N$ denote the kernel of the map $\T^2_A(F) \sra \T^2_A(M)$. Then $N$ is generated by
\[
\mas{e_i \ten f, \ f \ten e_i: \ \ 1 \leq i \leq k}.
\]
To show the linear independence of $\mas{m_i \ten m_j}_{i,j \leq k-1}$ we assume that we have an element $e \in N$ that is a linear combination of $\mas{e_i \ten e_j}_{i,j \leq k-1}$, and we need to show that $e=0$. We have
\[
e = \sum_{i=1}^k (a_i e_i \ten f + b_i f \ten e_i) = \sum_{i=1}^k \sum_{j=1}^k (a_i f_j + b_j f_i) e_i \ten e_j,
\]
and $a_i$ and $b_i$ satisfy the equations
\[
a_i f_k + b_k f_i = 0, \quad 1 \leq i \leq k.
\]
When $i=k$ we obtain $a_k f_k = - b_k f_k$ and so $b_k = -a_k$. Further, the fact that $\gcd(f_k,f_i) = 1$ for $1 \leq i \leq k-1$ gives us $f_i \mid a_i$ for all $i$. Thus $a_i = c_i f_i$ for all $i$, for some $c_i \in A$ and we thus have
\[
0 = a_i f_k + b_k f_i = a_i f_k - a_k f_i = c_i f_i f_k - c_k f_k f_i
\]
and so $c_i = c_k$ for all $i$. Hence
\[
e = c_k \sum_{i=1}^k (f_i e_i \ten f - f_i f \ten e_i) = c_k(f \ten f - f \ten f) = 0.
\]
Next, let $m \in \TS^2_A(M)$. We need to show that $m$ is generated by the elements (\ref{eq:TS2_UFD_generators1}). We may assume that $m$ is of the form
\[
m = \sum_{1 \leq i<j \leq k} a_{ij} m_i \ten m_j
\]
with $a_{ij} \in A$. Since $m \in \TS^2_A(M)$ it follows that
\[
\sum_{1 \leq i<j \leq k} a_{ij} (e_i \ten e_j - e_j \ten e_i) \in N
\]
and so
\[
\sum_{1 \leq i<j \leq k} a_{ij} (e_i \ten e_j - e_j \ten e_i) = \sum_{i=1}^k (x_i e_i \ten f + y_i f \ten e_i)
\]
where $x_i, y_i \in A$. We thus obtain the equalities $a_{ij} = x_i f_j + y_j f_i$ for $1 \leq i < j \leq k$. We can then write $m$ as
\[
m = \sum_{1 \leq i<j \leq k} a_{ij} m_i \ten m_j = \underbrace{\sum_{1 \leq i<j \leq k-1} a_{ij} m_i \ten m_j}_{ m'} + \sum_{i = 1}^{k-1} a_{ik} m_i \ten m_k = 
\]
\[
= m' + \sum_{i = 1}^{k-1} (x_i f_k + y_k f_i) m_i \ten m_k =
\]
\[
= m' + \sum_{i = 1}^{k-1} x_i m_i \ten (f_k m_k) + y_k \pa{\sum_{i = 1}^{k-1} f_i m_i} \ten m_k =
\]
\eqbeg
\label{eq:TS2_UFD}
= m' + \sum_{i = 1}^{k-1} x_i m_i \ten \pa{-\sum_{j=1}^{k-1} f_j m_j} + y_k (-f_k m_k) \ten m_k.
\eqend
Now the first two terms of (\ref{eq:TS2_UFD}) is in $L^{\fS_2}$, so these are linear combinations of the elements (\ref{eq:TS2_UFD_generators2}). This shows that $m$ is a linear combination of the elements (\ref{eq:TS2_UFD_generators1}). 
\end{proof}

\begin{para}[Determining surjectivity algorithmically]
Let $A$ be a ring and $M$ an $A$-module of finite presentation, given as the cokernel of a map $A^l \sra A^m$. Denote by $F$ the free module $A^m$, and let $e_1, \ldots, e_m$ be a basis for $F$. 

Then we can algorithmically determine whether the map $\TS^n_A(F) \sra \TS^n_A(M)$ is surjective. By Proposition \ref{prp:canmap_surj_inj} this is equivalent to the canonical morphism $\Ga^n_A(M) \sra \TS^n_A(M)$ being surjective.\\
\\
Consider the surjection $\pi: \T^n_A(F) \sra \T^n_A(M)$ and let $N \sse \T^n_A(F)$ be the kernel of $\pi$. Choose generators $\sig_1, \ldots, \sig_k$ for the symmetric group $\fS_n$, which we may view as $A$-module homomorphisms
\[
\sig_j : \T^n_A(F) \ra \T^n_A(F), \ \ j= 1, \ldots, k
\]
by 
\[
\sig_j(e_{i_1} \ten \ldots \ten e_{i_n}) = e_{i_{\sig_j^{-1}(1)}} \ten \ldots \ten e_{i_{\sig_j^{-1}(n)}}.
\]
For each homomorphism $\sig_j$ we construct the homomorphism $u_j = 1_{\T^n_A(F)} - \sig_j$, and we let $K_j = \ker u_j \sse \T^n_A(F)$.

We now have by definition 
\eqbeg
\label{eq:alg_def_TS}
\TS^n_A(F) = \bigcap_{j=1}^k K_j.
\eqend
Define maps $v_j: \T^n_A(F) \sra \T^n_A(M)$ for $j = 1, \ldots, k$ by the composition:
\[
\xymatrix{
\T^n_A(F) \ar[r]^{u_j} &  \T^n_A(F) \ar[r]^{\pi} & \T^n_A(M). 
}
\]
Let $L_j = \ker v_j$, and consider the intersection $L = \bigcap_{j=1}^k L_j$. Then we have that 
\[
\TS^n_A(M) = \pi(L) \sse \T^n_A(M).
\]
It is clear that we have an inclusion 
\[
\TS^n_A(F) + N  = \bigcap_{j=1}^k K_j + N \sse L
\]
and the question of the surjectivity of $\TS^n_A(F) \sra \TS^n_A(M)$ is now reduced to checking if $\pi(\TS^n_A(F))$ is strictly contained in $\pi(L)$. Finally we have that
\[
\pi(\TS^n_A(F)) = \pi(L) = \TS^n_A(M)
\]
if and only if
\eqbeg
\label{eq:criteria_surjectivity}
\TS^n_A(F) + N = L
\eqend
as submodules of the free module $\T^n_A(F)$. 

Suppose that the ring $A$ is a quotient ring of the form $A = R/I$ where $R$ is a polynomial ring in finitely many variables over $\QQ$ or $\ZZ/(p)$ for a prime $p \geq 2$, and $I$ is an ideal. Then the submodules $K_j$, $N$ and $L$ as well as the intersection (\ref{eq:alg_def_TS}) and the relation (\ref{eq:criteria_surjectivity}) can be explicitly calculated with computer algebra software such as Macaulay2.
\end{para}

\begin{exam}
\label{ex:canmap_not_surj_1}
Let $k$ be a field of characteristic $3$ and let $A = k[s,t]$ be the polynomial ring in two variables. Consider the free module $F = A^2$ with generators $e_1$, $e_2$ and the module $M = F/K$, where $K$ is the submodule generated by the element $s e_1 - te_2 \in F$.

We wish to show that the natural map $\TS^3_A(F) \sra \TS^3_A(M)$ is not surjective. By Proposition \ref{prp:canmap_surj_inj} this implies that the canonical morphism $\Ga^3_A(M) \sra \TS^3_A(M)$ is not surjective.\\
\\
Consider the element $u = s e_1 \ten e_1 \ten e_2 \in \T^3_A(F)$. Let $m_i$ denote the image of $e_i$ in $M$. We wish to show that the element $\bar{u} = s m_1 \ten m_1 \ten m_2$ is in $\TS^3_A(M)$. Since $s m_1 = t m_2$ in $M$, we have
\[
\bar{u} = s m_1 \ten m_1 \ten m_2 = t m_2 \ten m_1 \ten m_2 =  s m_2 \ten m_1 \ten m_1.
\]
This demonstrates that $\bar{u}$ is invariant under the action of $\fS_3$.\\
\\
Assume now that $\bar{u}$ is the image of an element $v \in \TS^3_A(F)$. Let $N$ denote the kernel of the projection map $\T^3_A(F) \sra \T^3_A(M)$. Then we have
\[
u = v + w
\]
where $w \in N$. Let $n = s e_1 - t e_2 \in F$. Then $N$ is generated by the elements
\[
\mas{n \ten e_i \ten e_j, \ e_i \ten n \ten e_j, \ e_i \ten e_j \ten n}_{i,j = 1,2}
\]
and so $w$ is a sum of the form
\[
w = \sum_{i,j} \pa{a_{ij} n \ten e_i \ten e_j + b_{ij} e_i \ten n \ten e_j
+ c_{ij} e_i \ten e_j \ten n},
\]
where $i,j = 1,2$ and $a_{ij}, b_{ij}, c_{ij} \in A$. Let $f: \T^3_A(F) \sra \T^3_A(F)$ be defined by $f = 1_{\T^3_A(F)} + \sig + \sig^2$, where $\sig$ is the homomorphism corresponding to the permutation $(1 \ 2 \ 3) \in \fS_3$. Then
\[
f(u) = f(v+w) = f(v) + f(w) = f(w)
\]
since $f(v) = 3v = 0$.

Also,
\eqbeg
\label{eq:fu}
f(u) = s (e_1 \ten e_1 \ten e_2 + e_1 \ten e_2 \ten e_1 + e_2 \ten e_1 \ten e_1)
\eqend
and
\eqbeg
\label{eq:fw}
f(w) = \sum_{i,j} d_{ij} \pa{ n \ten e_i \ten e_j + e_i \ten n \ten e_j + e_i \ten e_j \ten n},
\eqend
where $d_{ij} = a_{ij} + b_{ij} + c_{ij}$. Since $n = s e_1 - t e_2$, the coefficient in front of $e_1 \ten e_1 \ten e_2$ in (\ref{eq:fw}) is $-t d_{11} - s d_{12}$ and the coefficient in front of $e_1 \ten e_2 \ten e_2$ is $t d_{12} + s d_{22}$.

Comparing these coefficients with (\ref{eq:fu}) gives
\[
\arbeg{rcl}
-t d_{11} - s d_{12} & = & s,\\
t d_{12} + s d_{22} & = & 0.
\arend
\]
The first equation leads to $s(d_{12} + 1) = -t d_{11}$ and so $t \mid (d_{12} + 1)$. Thus $d_{12} = th - 1$ with $h \in A$. From the second equation we obtain $t d_{12} = - s d_{22}$ and hence $s \mid d_{12}$. This contradicts the fact that $d_{12} = th - 1$, and we thus conclude that $\bar{u}$ cannot be the image of an element of $\TS^3_A(F)$.
\end{exam}

\begin{rem}
\label{rem:canmap_not_surj_1}
It is possible to extend Example \ref{ex:canmap_not_surj_1} to characteristic $p \geq 3$ by making the following modifications: The field $k$ is of characteristic $p$, while $A = k[s,t]$ and $M = F/(s e_1 - t e_2)$ as before, where $F$ is a free module with basis $e_1, e_2$. The goal is now to show that the map $\Ga^p_A(M) \to \TS^p_A(M)$ is not surjective. We consider the element
\[
u = s e_1 \ten \ldots \ten e_1 \ten e_2 \in \T^p_A(F)
\]
and one shows that its image $\bar{u} \in \T^p_A(M)$ is in $\TS^p_A(M)$. To show that $\bar{u}$ is not the image of an element of $\TS^p_A(F)$ one follows the method in the example, replacing the permutation $\sig = (1 \ 2 \ 3)$ with the permutation $\sig = (1 \cdots p)$ and the function $f$ with the function
\[
f = 1_{\T^p_A(F)} + \sig + \sig^2 + \ldots + \sig^{p-1}.
\]
\end{rem}

\begin{exam}
\label{ex:canmap_not_surj_2}
Here we give an example of an $A$-module $M$ and a base extension $A \sra A'$ such that $\Ga^n_A(M) \sra \TS^n_A(M)$ is surjective but $\Ga^n_{A'}(M \ten_A A') \sra \TS^n_{A'}(M \ten_A A')$ is not surjective.

Let $k$ be a field of characteristic $2$ and let 
\[
A = k[x_1, x_2, x_3, y_1, y_2, y_3, z_1, z_2, z_3]
\]
be the polynomial ring in $9$ variables. Consider the free module $F=A^3$ with basis $e_1, e_2, e_3$ and let $n = z_1 e_1 + z_2 e_2 + z_3 e_3 \in F$. Define the $A$-module $M = F / \langle n \rangle$. Then the map $\Ga^2_A(M) \sra \TS^2_A(M)$ is surjective by Lemma \ref{lem:TS2_UFD}.

Consider the ideal $I \sse A$ generated by the coefficients of
\[
\sum_{i=1}^3 (x_i e_i \ten n + y_i n \ten e_i) - (x_1 z_2 + y_2 z_1) (e_1 \ten e_2 - e_2 \ten e_1) \in F \ten_A F.
\]
This ideal is then generated by the elements
\[
\mas{x_1 z_2 + y_2 z_1 + x_2 z_1 + y_1 z_2} \cup \mas{x_i z_j + y_j z_i: \ \mas{i,j} \ne \mas{1,2}}.
\]

Let $A' = A/I$ and let $M' = M \ten_A A'$ and $F' = F \ten_A A'$. We have a surjection $F' \sra M'$ mapping the basis $\mas{e_1,e_2,e_3}$ to a set of generators $\mas{m_1,m_2,m_3}$ of $M'$. 

We wish to show that $\Ga^2_{A'}(M') \sra \TS^2_{A'}(M')$ is not surjective, or equivalently by Proposition \ref{prp:canmap_surj_inj} that $\TS^2_{A'}(F') \sra \TS^2_{A'}(M')$ is not surjective.

Consider the element
\[
u = (x_1 z_2 + y_2 z_1) e_1 \ten e_2 \in \T^2_{A'}(F')
\]
and let $\bar{u} = (x_1 z_2 + y_2 z_1) m_1 \ten m_2 \in \T^2_{A'}(M')$ be the image.

By the construction of the ideal $I \sse A$ we have that
\[
(x_1 z_2 + y_2 z_1) (m_1 \ten m_2 - m_2 \ten m_1) = 0 \in \T^2_{A'}(M')
\]
and hence $\bar{u} \in \TS^2_{A'}(M')$. Our aim is to show that $\bar{u}$ is not the image of an element in $\TS^2_{A'}(F')$.

Suppose therefore that $u = v+w$ for some $v \in \TS^2_{A'}(F')$ and $w$ in the kernel of $\T^2_{A'}(F') \sra \T^2_{A'}(M')$. Thus
\[
w = \sum_{i=1}^3 (a_i e_i \ten n + b_i n \ten e_i)
\]
where $a_i, b_i \in A'$. Let $f: \T^2_{A'}(F') \sra \T^2_{A'}(F')$ be defined by $f = 1_{\T^2_{A'}(F')} + \sig$, where $\sig(e_i \ten e_j) = e_j \ten e_i$. Then applying $f$ to the equation $u = v+w$ and using the fact that $f(v) = 0$ we obtain
\[
f(u) = (x_1 z_2 + y_2 z_1) (e_1 \ten e_2 + e_2 \ten e_1) = \sum_{i=1}^3 c_i (e_i \ten n + n \ten e_i) = f(w),
\]
where $c_i = a_i + b_i$. This equation leads to three equations involving the coefficients $c_i$:
\eqbeg
\label{eq:3eqs_TS2_3gens}
\arbeg{rcl}
c_1 z_2 + c_2 z_1 & = & x_1 z_2 + y_2 z_1\\
c_1 z_3 + c_3 z_1 & = & 0\\
c_2 z_3 + c_3 z_2 & = & 0.
\arend
\eqend
We now introduce a multigrading of the polynomial ring $A$ by
\[
\mdeg(x_i) = \mdeg(y_i) = \mdeg(z_i) = (1,i), \quad i = 1,2,3.
\]
With respect to this multigrading the ideal $I \sse A$ is homogeneous, and so the grading carries over to the quotient ring $A'=A/I$.

Since the right hand side of (\ref{eq:3eqs_TS2_3gens}) is homogeneous of multidegree $(2,3)$ we have that these equations are satisfied when $c_i$ is replaced by its homogeneous part of multidegree $(1,i)$ for $i=1,2,3$.

Thus we may assume that
\[
c_i = \al_i x_i + \be_i y_i + \ga_i z_i, \quad \al_i, \be_i, \ga_i \in k, \quad i = 1,2,3.
\]
We will now show that the equations (\ref{eq:3eqs_TS2_3gens}) lead to a contradiction. When working in the ring $A'=A/I$ we will make the following reductions of binomials: 
\[
y_1 z_2 \sra x_1 z_2 + x_2 z_1 + y_2 z_1, \quad y_i z_j \sra x_j z_i, \ \mas{i,j} \ne \mas{1,2}. 
\]
Now consider an integer $i \in \mas{1,2}$. We work out the last two equations of (\ref{eq:3eqs_TS2_3gens}) as follows:
\[
c_i z_3 + c_3 z_i = \al_i x_i z_3 + \be_i y_i z_3 + \ga_i z_i z_3 + \al_3 x_3 z_i + \be_3 y_3 z_i + \ga_3 z_3 z_i =
\]
\[
= \al_i x_i z_3 + \be_i x_3 z_i + \ga_i z_i z_3 + \al_3 x_3 z_i + \be_3 x_i z_3 + \ga_3 z_3 z_i = 0.
\]
This gives
\[
\al_i = \be_3, \ \ \be_i = \al_3, \ \ \ga_i = \ga_3, \quad i \in \mas{1,2}.
\]
The first equation of (\ref{eq:3eqs_TS2_3gens}) now becomes
\[
c_1 z_2 + c_2 z_1 = \be_3 x_1 z_2 + \al_3 y_1 z_2 + \ga_3 z_1 z_2 + \be_3 x_2 z_1 + \al_3 y_2 z_1 + \ga_3 z_2 z_1 = 
\]
\[
= \be_3 x_1 z_2 + \al_3 (x_1 z_2 + x_2 z_1 + y_2 z_1) + \be_3 x_2 z_1 + \al_3 y_2 z_1 =
\]
\[
= (\al_3 + \be_3)(x_1 z_2 + x_2 z_1) \ne x_1 z_2 + y_2 z_1
\]
and this is the desired contradiction. The conclusion is that the element
\[
\bar{u} = (x_1 z_2 + y_2 z_1) m_1 \ten m_2 \in \TS^2_{A'}(M')
\]
is not the image of an element of $\TS^2_{A'}(F')$. Hence the canonical map $\Ga^2_{A'}(M') \sra \TS^2_{A'}(M')$ is not surjective.
\end{exam}

\begin{rem}
\label{rem:canmap_not_surj_2}
Example \ref{ex:canmap_not_surj_2} works by choosing the ideal $I \sse A$ to be the ideal defining the relation that the element $(x_1 z_2 + y_2 z_1) m_1 \ten m_2 \in \T^2_{A'}(M')$ is symmetric.

It might be possible to make a similar construction in characteristic $p \geq 2$ by choosing the ring $A$ to be a large polynomial ring and constructing the ideal $I \sse A$ to be the ideal defining the relation that an element of the form
\[
f m_1 \ten \ldots \ten m_1 \ten m_2 \in \T^p_{A'}(M')
\]
is symmetric, where $M = A^3/(z_1 e_1 + z_2 e_2 + z_3 e_3)$ as before, $A' = A/I$, and $f \in A$ is some polynomial. One might then be able to use methods similar to the ones in Example \ref{ex:canmap_not_surj_2} to show that $\Ga^p_{A'}(M') \to \TS^p_{A'}(M')$ is not surjective.

The map $\Ga^p_A(M) \to \TS^p_A(M)$ is probably surjective but to show this we would require a modification of Lemma \ref{lem:TS2_UFD} to deal with $n > 2$, and this we do not know how to do.
\end{rem}

\section{Symmetric tensors and base change}

In this section we give examples to show that the functor $\TS$ of symmetric tensors does not commute with base change in general.

\begin{defi}
Let $A \sra A'$ be a homomorphism of rings, and consider an $A$-module $M$. Denote by $M'$ the module $M \ten_A A'$ obtained by base extension to $A'$. We have a natural isomorphism
\[
\T_A^n(M) \ten_A A' \stackrel{\sim}{\ra} \T_{A'}^n(M')
\]
inducing a canonical map
\eqbeg
\label{eq:base_ext_map}
\TS_A^n(M) \ten_A A' \ra \TS_{A'}^n(M').
\eqend
\end{defi}

\begin{prop}
The base change morphism (\ref{eq:base_ext_map}) is an isomorphism in the following cases:
\begin{itemize}
\item[(i)] The element $n!$ is invertible in the ring $A$.
\item[(ii)] The $A$-module $M$ is flat.
\item[(iii)] The base extension $A \sra A'$ is flat.
\end{itemize}
\end{prop}

\begin{proof}
To show (i) and (ii) we consider the commutative diagram
\[
\xymatrix{
\Ga^n_A(M) \ten_A A' \ar[r] \ar[d] & \Ga^n_{A'}(M') \ar[d]\\
\TS^n_A(M) \ten_A A' \ar[r] & \TS^n_{A'}(M')
}
\]
where the top horizontal map is the map (\ref{par:Ga_base_change}), the bottom horizontal map is the base change morphism and the vertical maps are the canonical maps of Definition ~\ref{def:Ga_TS_canmap}. By (\ref{par:Ga_base_change}) the top horizontal map is an isomorphism and by Proposition ~\ref{prp:Ga_TS_iso} both vertical maps are isomorphisms. Hence the bottom horizontal map is an isomorphism.

To show (iii), let $\sig_1, \ldots, \sig_k$ be generators of $\fS_n$ regarded as morphisms $\sig_i: \T^n_A(M) \sra \T^n_A(M)$. Then $\TS^n_A(M)$ is the submodule of $\T^n_A(M)$ consisting of those $x \in \T^n_A(M)$ such that $\sig_i(x) = \sig_j(x)$ for all $i,j$. In other words, $\TS^n_A(M)$ is the inverse limit of the diagram
\[
\xymatrix{
\T^n_A(M) \ar@/^1pc/[rr]^{\sig_1} \ar@/_1pc/[rr]_{\sig_k} & \vdots & \T^n_A(M).
}
\]
The base extension functor $N \mapsto N \ten_A A'$ is exact since $A'$ is flat, and exact functors commute with finite inverse limits \cite[Def. 2.4.1]{sga4_grothendieck_prefaisceaux}. We have that $\T^n_A(M) \ten_A A' \cong \T^n_{A'}(M')$, so $\TS^n_{A'}(M')$ is the inverse limit of the diagram
\[
\xymatrix{
\T^n_A(M) \ten_A A' \ar@/^1pc/[rr]^{\sig_1 \ten 1} \ar@/_1pc/[rr]_{\sig_k \ten 1} & \vdots & \T^n_A(M) \ten_A A'.
}
\]
The fact that flat base extension commutes with finite inverse limits shows that the canonical map
\[
\TS^n_A(M) \ten_A A' \ra \TS^n_{A'}(M')
\]
is an isomorphism.
\end{proof}

\begin{exam}
\label{ex:TS_base_change_not_injective}
Here we give an example where the base change map is not \emph{injective}. Let the morphism of rings $A \sra A'$ and the $A$-module $M$ be as in Example \ref{ex:canmap_not_injective}. That is, $k$ is a field of characteristic $2$, the ring $A = k[s,t]$ is the polynomial ring and $M = A^2/ \langle s e_1 + t e_2 \rangle$ with $\mas{e_1, e_2}$ being the natural basis of $A^2$. Furthermore $A' = A[z]/(z(s+t))$.

We have the canonical commutative diagram
\[
\xymatrix{
\Ga^2_A(M) \ten_A A' \ar[r] \ar[d] & \Ga^2_{A'}(M') \ar[d]\\
\TS^2_A(M) \ten_A A' \ar[r] & \TS^2_{A'}(M')
}
\]
where the top horizontal map is an isomorphism by (\ref{par:Ga_base_change}) and the bottom horizontal map is the base change morphism. By Example \ref{ex:canmap_not_injective} the map $\Ga^2_A(M) \sra \TS^2_A(M)$ is injective, and since $M$ is generated by $2$ elements the map $\Ga^2_A(M) \sra \TS^2_A(M)$ is surjective as well by Lemma \ref{lem:2gensTS2}. Hence the leftmost map of the diagram is an isomorphism and in particular injective. 

However, by Example \ref{ex:canmap_not_injective} the rightmost vertical map is not injective. Thus the bottom horizontal map cannot be injective. Specifically, Example \ref{ex:canmap_not_injective} shows that the element 
\[
(m_1 \ten m_1 + m_2 \ten m_2) \ten zs \in \TS^2_A(M) \ten_A A'
\]
is nonzero and is mapped to zero in $\TS^2_{A'}(M')$.
\end{exam}

\begin{rem}
It might be possible to extend Example \ref{ex:TS_base_change_not_injective} to characteristic $p > 2$ by using Remark \ref{rem:canmap_not_injective}. With the notation of Remark \ref{rem:canmap_not_injective}, we have that $\Ga^p_{A'}(M') \to \TS^p_{A'}(M')$ is not injective. It is probably true that $\Ga^p_{A}(M) \to \TS^p_{A}(M)$ is an isomorphism, but this has not been proven.
\end{rem}

\begin{exam}
\label{ex:TS_base_change_not_surjective}
Here we give an example where the base change map is not \emph{surjective}. Let the morphism of rings $A \sra A'$ and the $A$-module $M$ be as in Example \ref{ex:canmap_not_surj_2}. That is, $k$ is a field of characteristic $2$, the ring $A = k[x_1,x_2,x_3,y_1,y_2,y_3,z_1,z_2,z_3]$ is the polynomial ring and $A' = A/I$ where $I$ is the ideal generated by the elements
\[
\mas{x_1 z_2 + y_2 z_1 + x_2 z_1 + y_1 z_2} \cup \mas{x_i z_j + y_j z_i: \ \mas{i,j} \ne \mas{1,2}}.
\]
The module $M$ is defined as $M = A^3/ \langle z_1 e_1 + z_2 e_2 + z_3 e_3 \rangle$, where $\mas{e_1,e_2,e_3}$ is the natural basis of $A^3$.

Then we have the canonical commutative diagram
\[
\xymatrix{
\Ga^2_A(M) \ten_A A' \ar[r] \ar[d] & \Ga^2_{A'}(M') \ar[d]\\
\TS^2_A(M) \ten_A A' \ar[r] & \TS^2_{A'}(M')
}
\]
where the top horizontal map is an isomorphism by (\ref{par:Ga_base_change}) and the bottom horizontal map is the base change morphism. By Example \ref{ex:canmap_not_surj_2} the leftmost vertical map is surjective while the rightmost vertical map is not surjective. Thus the bottom horizontal map cannot be surjective. Specifically, Example \ref{ex:canmap_not_surj_2} shows that the element
\[
(x_1 z_2 + y_2 z_1) m_1 \ten m_2 \in \TS^2_{A'}(M')
\]
is not the image of an element of $\TS^2_A(M) \ten_A A'$.
\end{exam}

\begin{rem}
It might be possible to extend Example \ref{ex:TS_base_change_not_surjective} to characteristic $p > 2$ by using Remark \ref{rem:canmap_not_surj_2}. With the notation of Remark \ref{rem:canmap_not_surj_2}, one may be able to show that $\Ga^p_{A'}(M') \to \TS^p_{A'}(M')$ is not surjective. The map $\Ga^p_{A}(M) \to \TS^p_{A}(M)$ is probably surjective, but this has not been shown.
\end{rem}

\section{From modules to algebras}
\label{sec:from_modules_to_algebras}

In the previous sections we have given examples of modules $M$ such that the canonical map is not an isomorphism and such that the symmetric tensors do not commute with base change. Here we extend the previous examples to algebras.

\begin{prop}
\label{prp: graded_ring_functors}
Let $A$ and $A'$ be rings and let $R = \bigoplus_{k \geq 0} R_k$ be a graded $A$-algebra. Let $F,G: \Mod{A} \sra \Mod{A'}$ be covariant functors and consider a natural transformation $\fii: F \sra G$. 

If the map $\fii_R : F(R) \sra G(R)$ is injective (resp. surjective), then the map $\fii_{R_k} : F(R_k) \sra G(R_k)$ is injective (resp. surjective), where $R_k$ denotes the $k$:th graded piece of $R$.
\end{prop}

\begin{proof}
We have a canonical inclusion map $R_k \sra R$ with a section $R \sra R_k$ given by the projection onto the $k$:th factor. Applying the functors $F$ and $G$ to the sequence $R_k \sra R \sra R_k$ gives a commutative diagram
\[
\xymatrix{
F(R_k) \ar[r] \ar[d]_{\fii_{R_k}} & F(R) \ar[r] \ar[d]_{\fii_R} & F(R_k) \ar[d]^{\fii_{R_k}}\\
G(R_k) \ar[r] & G(R) \ar[r] & G(R_k)
}
\]
of $A'$-modules. The composition of the left and right top horizontal arrow gives the identity, and likewise for the bottom horizontal arrows. Thus the left horizontal arrows are injective and the right are surjective.

Suppose that $\fii_R$ is injective. Then one concludes from the leftmost square that $\fii_{R_k}$ is injective. Next, if $\fii_R$ is surjective we conclude from the rightmost square that also $\fii_{R_k}$ is surjective.
\end{proof}

\begin{cor}
\label{cor:graded_ring_Ga_TS}
Let $A$ be a ring and $R = \bigoplus_{k \geq 0} R_k$ a graded $A$-algebra. Suppose that the canonical map $\Ga^n_A(R) \sra \TS^n_A(R)$ is injective (resp. surjective). Then the map $\Ga^n_A(R_k) \sra \TS^n_A(R_k)$ is injective (resp. surjective).
\end{cor}

\begin{proof}
In Proposition \ref{prp: graded_ring_functors} choose $A' = A$, $F = \Ga^n_A$ and $G = \TS^n_A$. Let $\fii: F \sra G$ be the canonical map.
\end{proof}

\begin{cor}
\label{cor:graded_ring_TS_base_change}
Let $A$ be a ring, $A \sra A'$ an $A$-algebra and $R = \bigoplus_{k \geq 0} R_k$ a graded $A$-algebra. Suppose that the canonical base change map $\TS^n_A(R) \ten_A A' \sra \TS^n_{A'}(R \ten_A A')$ is injective (resp. surjective). Then the map $\TS^n_A(R_k) \ten_A A' \sra \TS^n_{A'}(R_k \ten_A A')$ is injective (resp. surjective).
\end{cor}

\begin{proof}
In Proposition \ref{prp: graded_ring_functors} choose $F(\cdot) = \TS^n_A(\cdot) \ten_A A'$ and $G(\cdot) = \TS^n_{A'}(\ \cdot \ \ten_A A')$. Let $\fii: F \sra G$ be the base change morphism.
\end{proof}

Proposition \ref{prp: graded_ring_functors} and its corollaries extend the examples of the previous sections to algebras, by considering the symmetric algebra $\rmS_A(M)$ of an $A$-module $M$. 

\begin{exam}
\label{ex:modules_to_algebras}
Examples of rings $A$ and $A$-algebras $B$ and $A'$ such that
\begin{itemize}
\item[(a)] $\Ga^n_A(B) \sra \TS^n_A(B)$ is not injective. Let $k$ be a field of characteristic $p > 0$, and let $A = k[x]$ and $B = k$, where $A \sra B$ sends $x$ to $0$. Then by Example \ref{ex:rydh_canmap_not_injective} we have that $\Ga^p_A(B) \sra \TS^p_A(B)$ is not injective.
\item[(b)] $\Ga^n_A(B) \sra \TS^n_A(B)$ is not surjective. Here we choose $k$ a field of characteristic $p \geq 3$, $A = k[s,t]$ and $B = A[x,y]/(sx - ty)$. Then $B = \rmS_A(M)$, the symmetric algebra of the module $M$ of Example \ref{ex:canmap_not_surj_1}. Thus $\Ga^p_A(B) \sra \TS^p_A(B)$ is not surjective by Example \ref{ex:canmap_not_surj_1}, Remark \ref{rem:canmap_not_surj_1} and Corollary \ref{cor:graded_ring_Ga_TS}. 
\item[(c)] $\TS^n_A(B) \ten_A A' \sra \TS^n_{A'}(B \ten_A A')$ is not injective. Choose $k$ to be a field of characteristic $2$, $A = k[s,t]$ the polynomial ring, and $A' = A[z]/(z(s+t))$. Furthermore, let $B = A[x,y]/(sx + ty)$. Then $B = \rmS_A(M)$, the symmetric algebra of the module $M$ of Example \ref{ex:TS_base_change_not_injective}. Therefore $\TS^2_A(B) \ten_A A' \sra \TS^2_{A'}(B \ten_A A')$ is not injective by Example \ref{ex:TS_base_change_not_injective} and Corollary \ref{cor:graded_ring_TS_base_change}.
\item[(d)] $\TS^n_A(B) \ten_A A' \sra \TS^n_{A'}(B \ten_A A')$ is not surjective. Let $k$ be a field of characteristic $2$ and let $A = k[x_1, x_2, x_3, y_1, y_2, y_3, z_1, z_2, z_3]$, the polynomial ring in $9$ variables. Choose $A' = A/I$ where $I \sse A$ is the ideal of Example \ref{ex:TS_base_change_not_surjective}. Furthermore, let $B = A[u,v,w]/(z_1 u + z_2 v + z_3 w)$. Then $B = \rmS_A(M)$, the symmetric algebra of the module $M$ of Example \ref{ex:TS_base_change_not_surjective}. Therefore the base change map $\TS^2_A(B) \ten_A A' \sra \TS^2_{A'}(B \ten_A A')$ is not surjective by Example \ref{ex:TS_base_change_not_surjective} and Corollary \ref{cor:graded_ring_TS_base_change}.
\end{itemize}
\end{exam}

\begin{rem}
Consider the $A$-algebra $B$ of Example \ref{ex:modules_to_algebras} (b), with $p=3$.

It is not hard to show that the algebra $\T^3_A(B)$ is \emph{reduced}, and we thus have that $\TS^3_A(B)$ is reduced. It follows that the homomorphism
\eqbeg
\label{eq:red_not_surj}
\Ga_A^3(B)_{\mathrm{red}} \ra \TS_A^3(B)_{\mathrm{red}}
\eqend
is not surjective.

David Rydh \cite{rydh_representability_of_gamma_inprep} has shown that for any $A$-algebra $B$, the morphism
\[
\Ga_A^n(B)_{\mathrm{red}} \ra \TS_A^n(B)_{\mathrm{red}}
\]
is \emph{injective}, and also that the morphism
\eqbeg
\label{eq:specTS_specGa}
\spec(\TS_A^n(B)) \ra \spec (\Ga_A^n(B))
\eqend
is a universal homeomorphism with trivial residue field extensions. However, the example (\ref{eq:red_not_surj}) shows that despite this, we do not have an induced isomorphism on the reduced structures of the schemes $\spec(\TS_A^n(B))$ and $\spec(\Ga_A^n(B))$.
\end{rem}

\section{Algebra structures on exterior powers}
\label{sec:alg_struct_on_ext_pow}

In this section we discuss a problem related to the work of Laksov and Thorup in \cite{laksovthorup_determinantal}. Let $A$ be a ring and let $B = A[x]$ be the polynomial ring in one variable. In the article the authors consider a $\TS^n_A(B)$-module structure on the exterior product $\bigwedge^n_A(B)$ and use this to obtain formulas related to Schubert calculus for Grassmannians and the intersection theory of flag schemes. We give here an example to show that such a $\TS^n_A(B)$-module structure does not exist in general.

\begin{para}
Let $A$ be a ring and $B$ an $A$-algebra. Consider also a $B$-module $M$ viewed as an $A$-module by restriction of scalars. Recall that the exterior product $\bigwedge^n_A(M)$ is the $A$-module defined as the tensor product $\T^n_A(M)$ modulo the submodule generated by elements of the form $m_1 \ten \ldots \ten m_n$ with $m_i = m_j$ for some $0 \leq i < j \leq n$.

Note that $\T^n_A(B)$ has a structure of commutative $A$-algebra by the multiplication 
\[
(x_1 \ten \ldots \ten x_n) \cdot (y_1 \ten \ldots \ten y_n) = x_1 y_1 \ten \ldots \ten x_n y_n.
\]
The symmetric group $\fS_n$ acts on $\T^n_A(B)$ by $A$-algebra homomorphisms, and so $\TS^n_A(B) = \T^n_A(B)^{\fS_n}$ is a subalgebra of $\T^n_A(B)$. Moreover, the $A$-module $\T^n_A(M)$ is canonically a $\T^n_A(B)$-module by the rule
\[
(x_1 \ten \ldots \ten x_n) \cdot (m_1 \ten \ldots \ten m_n) = x_1 m_1 \ten \ldots \ten x_n m_n.
\]
We have a canonical surjection $\phi: \T^n_A(M) \sra \bigwedge^n_A(M)$ of $A$-modules and we ask for a $\TS^n_A(B)$-module structure on the exterior product $\bigwedge^n_A(M)$ such that the map $\phi$ is $\TS^n_A(B)$-linear.

Laksov and Thorup has shown that a unique such $\TS^n_A(B)$-module structure exists on $\bigwedge^n_A(M)$ when either $M$ or $B$ are free as $A$-modules, or 2 is invertible in $B$, see \cite[Prop. 1.3]{laksovthorup_determinantal}. Such a structure does not exist in general as shown by the example below.
\end{para}

\begin{lemma}
\label{lem:kernel_of_can_surjection}
Let $A$ be a ring and $B$ an $A$-algebra. Then the kernel of the map $\phi:\TS^2_{A}(B) \sra \bigwedge^2_{A}(B)$ is the image of the canonical morphism $\Ga^2_A(B) \sra \TS^2_A(B)$.
\end{lemma}

\begin{proof}
By the definition of $\bigwedge^2_{A}(B)$ we have that the kernel $K = \ker(\phi)$ is generated by all elements of the form $x \ten x$ with $x \in B$. By Proposition ~\ref{prp:canmap_surj_inj} the image $I = \im(\Ga^2_A(B) \sra \TS^2_A(B))$ is generated by all elements of the form $x \ten x$ and $x \ten x' + x' \ten x$ with $x,x' \in B$. Thus $K \sse I$ and the simple relation
\[
x \ten x' + x' \ten x = (x+x') \ten (x+x') - x \ten x - x' \ten x'
\]
shows that also $I \sse K$.
\end{proof}

\begin{exam}
Let $k$ be a field of characteristic 2 and let $A$ be the quotient ring $A = k[x_1, x_2, x_3, y_1, y_2, y_3, z_1, z_2, z_3]/I$ where $I$ is the ideal generated by the elements
\[
\mas{x_1 z_2 + y_2 z_1 + x_2 z_1 + y_1 z_2} \cup \mas{x_i z_j + y_j z_i: \ \mas{i,j} \ne \mas{1,2}}.
\]
This is the ring denoted by $A'$ in Example \ref{ex:canmap_not_surj_2}. Let $B$ be the $A$-algebra defined by $B = A[u,v,w]/(z_1 u + z_2 v + z_3 w)$. Then $B$ is a graded ring such that $B_1$ is the $A$-module denoted by $M'$ in Example \ref{ex:canmap_not_surj_2}. The canonical map
\[
\Ga^2_{A}(B) \sra \TS^2_{A}(B)
\]
is therefore not surjective by Example \ref{ex:canmap_not_surj_2} and Corollary \ref{cor:graded_ring_Ga_TS}. Thus by Lemma \ref{lem:kernel_of_can_surjection} there is an element $\eta \in \TS^2_{A}(B)$ that does not map to zero via the canonical map $\phi:\T^2_{A}(B) \sra \bigwedge^2_{A}(B)$. Suppose there is a $\TS^2_{A}(B)$-module structure on $\bigwedge^2_{A}(B)$ making $\phi$ into a $\TS^2_{A}(B)$-module homomorphism. Then
\[
0 = \eta \cdot \phi(1 \ten 1) = \phi(\eta \cdot (1 \ten 1)) = \phi(\eta) \ne 0
\]
which is our desired contradiction.
\end{exam}


\bibliographystyle{dary}

\bibliography{papers}

\end{document}